\documentclass[12pt]{amsart}
\usepackage{fullpage}
\usepackage{amssymb}
\usepackage{amsmath}
\usepackage{fullpage}
\usepackage{array}      
\usepackage{amsmath,amsthm,amscd,amsfonts,amssymb}
\usepackage{latexsym} \theoremstyle{plain} \newtheorem{exam}{Example}
\newtheorem{defi}[exam]{Definition} \newtheorem{thm}[exam]{Theorem}
\newtheorem{lem}[exam]{Lemma} \newtheorem{prop}[exam]{Proposition}
\newtheorem{rem}[exam]{Remark} 
 \newcommand{\half}{{\frac{1}{2}}}
\newcommand{\pf}{{\noindent \bf Proof: }}
\newcommand{\enpf}{\begin{flushright} $\Box$ \end{flushright}}

\newtheorem{lem:peter}{Lemma}

\newcommand\n{{\bf n}}
\newcommand\m{{\bf m}}
\newcommand\x{{\bf x}}
\newcommand{\Schr}{Schr\"{o}dinger}
\newcommand{\rf}[1]{(\ref{#1})}
\newcommand{\ds}{\displaystyle}
\newcommand{\beq}{\begin{equation}}
\newcommand{\eeq}{\end{equation}}
\newcommand{\ba}{\begin{array}}
\newcommand{\ea}{\end{array}}
\newcommand{\Sup}{\displaystyle \sup}
\newcommand{\Sum}{\displaystyle \sum }
\let\cal=\mathcal
\begin{document}
\input{amssym.def}
\newcommand{\ctext}[1]{\makebox(0,0){#1}}
\newcommand{\beqa}{\begin{eqnarray}}
\newcommand{\eeqa}{\end{eqnarray}}
\newcommand{\beqal}{\begin{eqnarray}\label}
\newcommand{\ca}{\chi_\alpha}
\newcommand{\mb}{\mathbb}
\newcommand{\om}{\Omega^1_X}
                                          
\newcommand{\cod}{ codim(X \backslash U) \ge 2}
\newcommand{\bc}{\begin{center}} 
\newcommand{\ec}{\end{center}}

\def\hpuk{$\mathrm H ^0(U,E_{P}(\frak u)\otimes \Omega^1)$}
\def\hppk{$\mathrm H ^0(U,E_{P}(\frak p)\otimes \Omega^1)$}
\def\hplk{$\mathrm H ^0(U,E_{P}(\frak l)\otimes \Omega^1)$}
\def\hapuk{$\mathrm H ^1(U,E_{P}(\frak u)\otimes \Omega^1)$}
\def\hpauk{$\mathrm H ^0(U,E_{P_1}(\frak u)\otimes \Omega^1)$}
\def\hpapak{$\mathrm H ^0(U,E_{P_1}(\frak p_{1})\otimes \Omega^1)$}
\def\hpaqk{$\mathrm H ^0(U,E_{P_1}(\frak q)\otimes \Omega^1)$}
\def\hapauk{$\mathrm H ^1(U,E_{P_1}(\frak u)\otimes \Omega^1)$}
\def\ra{$\rightarrow$}
\def\raf{$\stackrel{f}{\rightarrow}$}
\def\rafa{$\stackrel{f_1}{\rightarrow}$}   
\def\rag{$\stackrel{g}{\rightarrow}$}   
\def\raga{$\stackrel{g_1}{\rightarrow}$}
\def\dai{$\stackrel{\downarrow}{i}$}
\def\uda{$\Vert$}
\def\daj{$\stackrel{\downarrow}{j}$}
\def\Et{\tilde{E^{\gamma}_P}}
\def\Gt{$\tilde{\gamma}$}
\newcommand{\da}{\downarrow}    
\def\a{$\mathcal O_X$}
\def\b{$E_{G}(\frak g) \otimes \Omega^1_U$}
\def\c{$E_{P}(\frak p) \otimes E_{P}(\frak p)^* \otimes \Omega^1_U$}
\def\d{$E_{G}(\frak g) \otimes E_{G}(\frak g)^* \otimes  \Omega^1_U$}
\newcommand{\aeg}{E_G(\frak g)}
\newcommand{\Lg}{\frak g}
\newcommand{\aep}{E_P(\frak p)}
\newcommand{\T}{\sigma^*T_{E/P}}
\newcommand{\aegk}{E_G(\frak g) \otimes \Omega^1}
\newcommand{\lr}{\longrightarrow}
\newcommand{\Ll}{\longleftarrow}
\newcommand{\Wi}{\mathcal W_i}
\newcommand{\Wj}{\mathcal W_j}
\newcommand{\adt}{ad(\theta)}
\newcommand{\mr}{\frak}
\newcommand{\mc}{\mathcal}   
\setlength{\unitlength}{0.1mm}     
\newcommand{\bi}{\begin{itemize}}
\newcommand{\ei}{\end{itemize}}  
\def\hpuk{$\mathrm H ^0(U,E_{P}(\frak u)\otimes \Omega^1_{U'})$}
\def\hppk{$\mathrm H ^0(U,E_{P}(\frak p)\otimes \Omega^1_{U'})$}
\def\hplk{$\mathrm H ^0(U,E_{P}(\frak l)\otimes \Omega^1_{U'})$}
\def\hapuk{$\mathrm H ^1(U,E_{P}(\frak u)\otimes \Omega^1_{U'})$}
\def\hpauk{$\mathrm H ^0(U,E_{P_1}(\frak u)\otimes \Omega^1_{U'})$}
\def\hpapak{$\mathrm H ^0(U,E_{P_1}(\frak p_{1})\otimes \Omega^1_{U'})$}
\def\hpaqk{$\mathrm H ^0(U,E_{P_1}(\frak q)\otimes \Omega^1_{U'})$}
\def\hapauk{$\mathrm H ^1(U,E_{P_1}(\frak u)\otimes \Omega^1_{U'})$}
\def\ra{$\rightarrow$}
\def\raf{$\stackrel{f}{\rightarrow}$}
\def\rafa{$\stackrel{f_1}{\rightarrow}$}
\def\rag{$\stackrel{g}{\rightarrow}$}
\def\raga{$\stackrel{g_1}{\rightarrow}$}
\def\dai{$\stackrel{\downarrow}{i}$}
\def\uda{$\Vert$}                                 
\def\daj{$\stackrel{\downarrow}{j}$}
\title{On Harder-Narasimhan reductions for Higgs principal bundles}

\author{Arijit Dey}
\address{Institute of Mathematical Sciences}
\email{arijit@imsc.res.in}
\author{R.Parthasarathi}
\address{Chennai Mathematical Institute}
\email{partha@cmi.ac.in}

\maketitle

\begin{abstract}
The existence and uniqueness of H-N reduction for the 
Higgs principal bundles over nonsingular projective variety 
is shown in this article. We also extend the notion of H-N
reduction for $(\Gamma, G)$-bundles and ramified $G$-bundles over a smooth
curve.
\end{abstract}
\section{Introduction}

Let $X$ be any smooth projective variety over an algebraic closed field $k$ of
characteristic zero and $G$ be any reductive algebraic group over $k$. The 
problem of the existence and uniqueness of 
Harder-Narasimhan reductions (henceforth called briefly as 
H-N reduction)
for principal $G$-bundles 
was solved by Atiyah and Bott in \cite {AtB}. The notion of semistability that 
was used in this paper was using the adjoint representaion, 
i.e $E$ is semistable iff $E(\frak g)$ is semistable as a vector bundle.

In an earlier article by A.Ramanathan (\cite {AR2}) this question of 
H-N reduction was posed as a problem intrinsically on $E$ 
for its applicability to the positive characteristic cases as well.  
The problem in this setting was solved by Behrend \cite {BR} and in 
positive characteristics Behrend had a conjecture 
which was verified by Mehta and Subramaniam. 
Biswas and Holla gave a different approach to the solution following 
essentially the broad setting in Ramanathan's paper in \cite {BH}.

The aim of this paper is to generalise the methods of Biswas and Holla 
to give an unified approach to the case of principal bundles 
with Higgs structure on smooth projective varieties as well 
as the case of ramified bundles on smooth curves 
(as defined in \cite{BBN})(see \S 2,\S 5. for definitions).
Recall that for the case of Higgs vector bundles on smooth 
projective varieties the existence and uniqueness of H-N filtrations 
is
proved in \S 3 of Simpson \cite{S2}. The case of ramified bundles 
similarly generalise the existence and uniqueness of 
H-N filtrations for parabolic vector bundles (proved in \cite{MS}).
 
By a {\it principal object} in this 
introduction we mean a principal Higgs bundle on smooth projective varieties
or a ramified bundle on smooth projective curves.(See \S 1 for the
definition of H-N reduction's)

The main theorem of this paper is:
{\thm Let $E$ be a principal $G$ object on $X$. Then there exists a canonical
H-N reduction $(P,\sigma)$ where $P$ is a parabolic subgroup of $G$ 
and 
$\sigma : X \lr E/P$ is a section of the associated 
fibre bundle $E/P$ over $X$. The H-N reduction is unique in natural sense.
}

The brief layout of this paper is as follows; in \S 2 we define the
notions which are necessary for our problem,in \S 3 we prove the 
existence of the H-N reduction for Higgs principal bundles, in \S 4
we prove that for any finite dimensional representation $W$ of $G$ $E(W)$
is Higgs semistable iff $E$ is Higgs semistable, in \S 5 we prove
the uniqueness of the H-N reduction for Higgs principal bundles, in 
\S 6
we prove H-N reduction for ramified $G$-bundles.

This problem was posed by V.Balaji, we are very grateful to 
V.Balaji and
D.S.Nagaraj for constant encouragement
throughout this work. This paper owes a lot to the stimulating discussions
we had with them. Finally we would like to thank the referee for 
going through this paper very carefully and providing many suggestions.

\section{Preliminaries}
 Throughout this paper, unless
otherwise
stated, we have the following
notations :
\begin{enumerate}
\item{}
 We work over an algebraically closed field $k$ of characteristic $0$. 
\item{}  $G$ will stand for a connected reductive algebraic group (i.e,
for the 
linear algebraic group $G$,
unipotent radical of $G$, denoted by $R_u(G)$, is trivial). 
\item{}  $X$ is an irreducible smooth projective variety over k of
dimension
$d$ $\ge 2$. 
\item{} Let $X$ be embedded in $\mb P^N$ for some positive integer N, 
which is equivalent
to
fixing a very ample line bundle on $X$, we call it $H$ which we fix it
throughout the paper.
\end{enumerate} 
Let $E$ be a torsion-free coherent sheaf of rank $r$ over $X$.
Then the degree of $E$ is defined as $\deg(E) = c_1(E) \cdot H^{d-1}$, where
$c_1(E) \in \mathrm H^2(X,\mb Z)$ is the first Chern class.
Let $U$ be an open subscheme containing all points of codimension one,
$(\cod)$ such that $E|_U$ is locally free. Then ${\wedge}^r(E\mid_U)$
is an invertible sheaf on $U$, corresponding to a divisor class $D$ on 
$U$. Then we have $c_1(E\mid_U) = c_1({\wedge}^r E\mid_U) = D$.
From the functoriality of Chern classes, $c_1(E) = c_1(E\mid_U) = D$.

Let $U$ be an open subset of $X$ such that $\cod$, then $i_* \mathcal O_U 
\cong \mathcal O_X$, further if $\mathcal F$ is a locally free sheaf on $U$
then $i_*\mathcal F$ is a reflexive sheaf over $X$. 

It follows from the above that if two bundles agree at all points of
codimension one, they have the same degree. Also, we can talk about the
degree of torsion-free sheaf on $U$(since $\cod$). Let $E$ be a 
torsion-free
sheaf on $U$, then we define $\deg(E) = \deg(i_*(E))$, where $i$ is the
inclusion
$ i: U \longrightarrow X $. Further, if
$0\rightarrow F\rightarrow E\rightarrow G\rightarrow 0$ is an
exact sequence of torsion-free sheaves on $X$, then 
$\deg(E) = \deg(F) + \deg(G)$, also $\deg E^* = -\deg(E)$.
Again from the functoriality of Chern classes, for exact sequence of
locally free sheaves over $U$ we have the same.
For if, $0\rightarrow F\rightarrow E\rightarrow G\rightarrow 0$ is an
exact sequence of locally free sheaves on $U$, then we have $\deg(i_*E) =
\deg(i_*F) + \deg(i_*G)$.Therefore it follows that, 
$\deg(E) = \deg(F) + \deg(G)$
on $U$. Now we can define for a torsion-free sheaf $E$ on $X$, 
\[
\mu(E) = \frac{\deg(E)}{rank(E)} \in \mb Q
\]
A torsion-free sheaf $E$ is said to be $\mu$-semistable if for all 
coherent subsheaves $F \subseteq E$ we have  
\begin{center}
$ \mu(F) = \frac{\deg(F)}{rank(F)}\le \frac{\deg(E)}{rank(E)} = \mu(E)$.
\end{center}
Note that it is enough to check only for subsheaves whose quotients
are torsion-free.
Let $E$ be a torsion-free coherent sheaf on $X$
and $\theta: E\rightarrow E\otimes {\Omega^1_X}$, 
an ${\mathcal O}_X$--linear homomorphism of sheaf of modules, such that
$\theta\wedge\theta = 0$ in $E\otimes \Omega^2_X$, then the
pair $(E,\theta)$ is called \emph{Higgs
sheaf on $X$}. If $E$ is locally free, then the pair
$(E,\theta)$  is called \emph{Higgs bundle on $X$}. Let $F$ be a subsheaf
of $E$ such that $\theta\mid_F:F\rightarrow F\otimes \om$, then we say
$(F,\theta\mid_F)$ is a Higgs subsheaf of $(E,\theta)$. Further,
$(E,\theta)$ is called \emph{Higgs semistable}
if for every Higgs subsheaf $F\subseteq E$ we have, $\mu(F)\le\mu(E)$
. We also have the notion of Chern
classes of the principal $G$
bundle $E$ over $X$; let $c_i(E) \in H^*(X, \mb Z)$ be the 
Chern classes of the principal bundle
$E$. We will say that $c_i(E) = 0$ for any $i$ if $c_i(E(W)) = 0$ for
every $W \in Rep(G)$, where $Rep(G)$ is the category 
of finite dimensional representations of $G$ (\cite{BB} Rem.2.3,).

The aim of this paper is to find Higgs compatible H-N reduction
for Higgs principal $G$ bundle $E$ on $X$.  
Now we recall some facts which we need in our work.   
\begin{enumerate}
\item{} Let $E$ be a principal $G$-bundle on $X$. Let $Y$ be any
quasi-projective variety on which $G$ acts from the left, then we
define by $E(Y)$ as the associated fiber bundle with fiber type $Y$ which
is the following object: $E(Y) = (E\times Y)/G$ for the twisted action
of $G$ on $(E\times Y)$ given by $g.(e,y) = (e.g^{-1},g.y)$ where
$g,e,y$ are in $G,E,Y$ respectively.
 
\item{} Let $H$ be a closed subgroup of $G$. A principal $G$ bundle
$E$ on $X$ is said to have a $H$ structure or equivalently a reduction of 
structure group to $H$ if there exists an open subset $U$ of $X$ 
which contains all the points of codimension one $\cod$ and a
section $\sigma:  U \longrightarrow E(G/H)$. 
In other words a reduction of structure group to $H$ means giving a 
principal $H$ bundle $E_H$ over an open set $U$ which satisfies  $\cod$
such that $E_H(G) \cong E\mid_U$.
Given a reduction of structure group over an open subset $U$ with $\cod$
there is a unique maximal open set over which the reduction extends.
\item{} Let $H$ be a closed subgroup of $G$ and $U$ be an open subset of $X$.
 A principal $G$ bundle
$E$ on $U$ is said to have a $H$ structure or equivalently a reduction of
structure group to $H$ if there exists an open subset $U'$ of $U$
which contains all the points of codimension one $ codim(U \backslash U') 
\ge 2$ and a
reduction $\sigma:  U' \longrightarrow E(G/H)$.
In other words, a reduction of structure group to $H$ means giving a
principal $H$ bundle
$E_H$ on $U'$ such that $E_H(G) \cong E\mid_{U'}$.           
\item{} Let $\pi:E\rightarrow X$ be a principal $G$ bundle, then $E$
is called a semistable bundle if for every parabolic subgroup $P$ of $G$ and
any reduction of structure group $(P,\sigma)$ where $\sigma:
U\longrightarrow E(G/P)$ is a reduction with $\cod$
and for every dominant character $\chi$ of $P$ the associated line
bundle  $i_{*}L_{\chi}$ on $X$
 has degree $\le 0$ .      
\item{} Let $E$ be a principal $G$ bundle on $X$ and $\theta \in
\mathrm H^0(X,E(\frak {g})
\otimes \om)$,( $\frak {g}$ is the Lie algebra of $G$) with
$\theta\wedge\theta = 0$ in $E(\frak {g})\otimes \Omega^2_X$. Then the
pair $(E,\theta)$ is called a Higgs principal $G$ bundle on $X$.
\end{enumerate}
{\defi 
\label{Hred} Let $(E,\theta)$ be a Higgs principal 
$G$ bundle 
on $X$. Let
$(H,\sigma)$ be a reduction of $E$ to $H$ over an open subset $U$ of $X$
with
$\cod$ where
$H$ is a connected closed subgroup of
$G$  and $\theta_H \in\mathrm H^0(U, E_H(\frak {h}) \otimes \Omega^1_U$)
such that the following  diagram
commutes.
\begin{center}
\begin{tabular}{ccc}
${\mathcal{O}}_U$&$\stackrel{\theta}{\longrightarrow}$&  $E({\frak {g}}) 
\otimes \Omega^1_U$ \\
&$\stackrel{\theta_{H}}{\searrow}$&$\uparrow$ \\
&&$E_{H}({\frak {h}}) \otimes \Omega^1_U$ \\
&&$\uparrow$ \\
&&$0$ 
\end{tabular} 
\end{center}
Then the quadruple $(H,\sigma,\theta_H,U)$ is called Higgs reduction of
$E$ to $H$} .
{\rem  Let $(E,\theta)$ be a Higgs principal $G$ bundle over $U$, where
$U$ is an open subset of $X$ with $\cod$ . We can define
the Higgs reduction of $E$ to a subgroup $H$ in the same way as above. 
}
{\defi Let $(E,\theta)$ be a Higgs principal $G$ bundle over $X$, it is
 said to be Higgs semistable if  
for any Higgs reduction $(P,\sigma, \theta_P,U_P)$ to a  parabolic   
and for any dominant
character $\chi$ of $P$ the line bundle $L_\chi$ on $U_P$ has non-positive
degree which is same as saying the degree of the reflexive sheaf $i_*L_\chi$ 
over $X$
has non positive degree.
} 
{\defi  Let $(E,\theta)$ be a Higgs principal $G$ bundle over $U$, where $U$
is an open subset with $\cod$. Then it is
 said to be Higgs semistable if  
for any Higgs reduction $(P,\sigma, \theta_P,U_P)$ to a  parabolic, where 
 $ codim(U \backslash U_P) \ge 2$   
and for any dominant
character $\chi$ of $P$ the line bundle $L_\chi$ on $U_P$ has non-positive
degree which is same as saying the degree of the reflexive sheaf $i_*L_\chi$ 
over $X$
has non positive degree.

}
{\lem
\label{para} Let $(E,\theta)$ be a Higgs principal $G$ bundle over $X$.
Then  $(E,\theta)$ is Higgs semistable if and only if for any maximal parabolic
subgroup  $P$ of $G$ and for any  Higgs reduction $(P,\sigma,
\theta_P,U_P)$ 
we have $\deg i_* \sigma^*(T_{G/P}) \ge 0$ where $T_{G/P}$ is the relative 
tangent bundle for the projection $E\mid_U\rightarrow E\mid_U(G/P)$ and
$i$ is the inclusion $U \hookrightarrow X$.
}

{\pf Let $P$ be a maximal parabolic
subgroup of $G$. Then consider the following exact sequence of 
vector bundles on $U_P$; $0\rightarrow E_P(\frak p)\rightarrow  E_P(\frak
g)\rightarrow E_P(\frak
g/\frak p)\rightarrow 0$  given by the exact sequence 
$0\rightarrow \frak p\rightarrow
\frak g\rightarrow \frak g/\frak p\rightarrow 0$. Since $G$ is reductive 
there is a non-degenerate bilinear form
on $\frak g$ invariant under $G$; therefore we see that  $\deg E_P(\frak g) = 0$.
Therefore, $\deg E_P(\frak p)\le 0$ iff $\deg E_P(\frak g/\frak p)\ge 0$. 
So the 'only if' part is trivial. The other way can be proved
following the
proof of lemma 2.1. in \cite{AR1}

{\lem
\label{G&V}
When $G$ is $GL(n,k)$, the Higgs semistability of Higgs principal  
$G$
bundle $E$ on $X$ 
is equivalent to the Higgs semistability of the associated Higgs 
vector bundle
$E(k^n)$
by the standard representation.
}

{\pf 
Let $P$ be a maximal parabolic subgroup of $G$. Let
$(P,\sigma_P,\theta_P,U_P)$ be a Higgs reduction of $E$ to $P$, where
$P$ is parabolic and let
$E_P$ be the corresponding $P$ bundle. Let us consider the flag
$0\subset k^r\subset k^n$ corresponding to the parabolic subgroup $P$. Let
$V=E_P(k^n),W = E_P(k^r)$; then $W$ is a Higgs subbundle of $V$ on $U_P$
(with Higgs
structure induced by $\theta_P$). Note that $E_P(\frak g/\frak p)
= V/W\otimes W^*$ and Higgs  structures on both sides are same.
(since they are coming from $\theta$ and $\theta_P$).
Also $\mu(W) \le \mu(V)$ is equivalent to $\mu(W) \le
\mu(V/W)$ on $U_P$.\\
\[
 \deg((V/W)\otimes W^*) = -\deg W . rank(V/W) + rankW . \deg(V/W) 
\]
\[ 
\hspace{4cm}           = (\mu(V/W)-\mu(W)) . rank(V/W) .
rank(W)  
\]
From the above equation it is easy to see that $\mu(W) \le \mu(V)$ is
equivalent to $\deg((V/W)\otimes W^*) \ge 0$. Since $E_P(\frak g/\frak p)
= (V/W)\otimes W^*$ we conclude that $\deg(E_P(\frak g/\frak p)) \ge 0$ iff
$\mu(W) \le \mu(V)$.

Suppose $E$ is Higgs semistable. then we have, $\deg {i_P}_*(E_P(\frak
g/\frak p)) \ge 0$ on $X$ for any Higgs reduction of $(P,\sigma_P,\theta_P,U_P
)$ 
to a  maximal parabolic subgroup $P$ where $i_P$ is the inclusion of
$U_P$ to $X$. Let $W$ be a Higgs torsion-free sheaf of $V$ on $X$. We can
choose an open subset $U$ with $\cod$ such that $W$ is a Higgs subbundle
of
$V$ on $U$. Then $W = E_P(k^r)$ for some Higgs reduction
$(P,\sigma_P,\theta_P)$ corresponding to a flag $0 \subset  k^r\subset
k^n$ on $U$. Then $\deg{i_P}_*(E_P(\frak g/\frak p)) \ge 0$ implies that
$\mu(W) \le \mu(V)$ 
on $U$. Thus $V$ is Higgs semistable. Converse is also true.


Given a torsion-free coherent sheaf $E$ over $X$, there is  a unique
filtration of subsheaves
\begin{center}
$0=E_0\subset E_1\subset E_2\subset \cdots\cdots\cdots\subset E_{l-1}
\subset E_l=E$
\end{center}
which is characterised by the two conditions that all
the sheaves $\frac{E_i}{E_{i-1}}, i\in[1,l]$, are semistable torsion-free sheaves and
the $\mu\big(\frac{E_i}{E_{i-1}}\big)$'s
are
strictly decreasing as $i$ increases. This filtration is known as
\emph{Harder-Narasimhan filtration} for $E$
The following is the definition of H-N reduction for principal bundles in 
the literature(cf. \cite{AAB}). 
{\defi Let $E$ be a principal $G$ bundle on $X$ and $(P,\sigma_P,U_P)$ be
a reduction of structure group of $E$ to a parabolic subgroup $P$
of $G$, then this reduction is called H-N reduction 
if the following two conditions hold: 
\begin{enumerate}
\item{}If $L$ is the Levi factor of $P$ then the principal $L$ bundle
 $E_P\;{\times_P} \;L$ over $U_P$ is a semistable $L$ bundle.
\item{} For any dominant character $\chi$ of $P$ with respect to some
  Borel subgroup $B \subsetneq P$ of $G$, the associated line bundle
  $L_\chi$ over $U_P$ has degree $> 0$.
\end{enumerate}
}                                      
For $G=GL(n,k)$ a reduction $E_P$ gives a filtration of the rank $n$
vector bundle associated to the standard representation. It is easy to see
that $E_P$ is canonical in the above sense iff the corresponding
filtration of the associated vector bundle filtration coincides with its
Harder Narasimhan filtration.
{\lem
Let $\rho: G\longrightarrow H$ be the homomorphism of connected
reductive algebraic groups and let $(E,\theta)$ be the Higgs $G$ bundle
on $X$.
Then the associated bundle $E_H : = E(H)$ gets a Higgs structure.
}
 
{\pf The representation $\rho $ induces a morphism $\rho': E(\frak
{g})\longrightarrow    
 E(\frak {h})$ of vector bundles on $X$.  
where ${\frak {g}}$ and ${\frak {h}}$ are Lie algebras of $G$ and $H$
respectively. So we define Higgs structure on $E(H)$ denoted by
$\theta_H :=(\rho'\otimes{id})\circ\theta$
where $\rho'\otimes{id} : E(\frak {g})\otimes{\om}\longrightarrow
E(\frak {h})\otimes{\om}$ and $\theta:  \mathcal{O}_X\longrightarrow
E(\frak {g})\otimes{\om}$.

Now we give a definition of a Higgs 
compatible  H-N reduction for Higgs principal bundles.   
{\defi Let $(E,\theta)$ be a Higgs principal $G$ bundle on $X$. Then a
Higgs
reduction $(E_P,\sigma_P,\theta_P,U_P)$ (Definition \ref{Hred}) 
is called a
Higgs compatible H-N
reduction if the following conditions hold:
\begin{enumerate}
\item{}The Higgs bundle $E_L$  is Higgs semistable on $U_P$, where $L$ is
the Levi factor of
$P$.
\item{}For all dominant character $\chi$ of $P$ with respect to some
Borel subgroup $B\subset P\subset G$ the associated line bundle
$L_\chi$
has positive degree on $U_P$.
\end{enumerate}
}

From now onwards we call a Higgs compatible H-N reduction as a
Higgs H-N reduction. 
 Suppose $H_1$, $H_2$
are closed connected subgroups
of $G$. Let $\sigma_1$, $\sigma_2$ be reductions of structure group of $E$
to
$H_1$ and $H_2$ respectively. To this data we can associate subgroup
schemes
$\sigma _{i}^*E(H_i)$ of $E(G)$ and their Lie algebras $\sigma_i^*E(\frak
h_i)$,
which are sub-bundles of $E(\frak g)$. For details we refer 
to \cite {BH}.

We recall some facts which we use in our proof of Higgs H-N
reduction for Higgs principal bundles:

\begin{enumerate}
\item{} Let $(E_1,\theta_1)$ and $(E_2,\theta_2)$
  be two Higgs semistable vector bundles on smooth projective variety
  $X$. Then their tensor product $(E_1\otimes E_2,\theta_1 \otimes id
  + \theta_2 \otimes id)$ is also Higgs semistable.(cf.\cite{S1})
  
\item{} For two torsion-free torsion-free coherent sheaves $E_1$ and
  $E_2$ over $X$, the equality\\ $\mu_{min}(E_1 \otimes E_2) =
  \mu_{min}(E_1) + \mu_{min}(E_2)$\\ is valid. Similarly, we have
  $\mu_{max}(E_1 \otimes E_2) = \mu_{max}(E_1) +
  \mu_{max}(E_2)$.(cf.\cite{AB})
\item{} Let $E_1$ and $E_2$ be two torsion-free coherent
sheaves over $X$.(cf.\cite{AB})                   
\begin{itemize}
\item If $\mu_{min}(E_1) > \mu_{max}(E_2)$, then $H^{0}(X,Hom_{\mathcal
O_{X}}(E_1,E_2)) = 0 $. 
\item If there is a surjective $\mathcal O_X$ linear homomorphism
   $\phi:E_1 \rightarrow E_2$ , 
 then  $\mu_{min}(E_1) \le \mu_{min}(E_2)$.

\end{itemize}
\end{enumerate}
}
\section{Proof of existence of a Higgs H-N reduction}

For a Higgs $G$ bundle $(E,\theta)$ on $X$ (following \cite {BH}) we 
define an 
integer 
$d_E$ as follows:

\begin{center}
$d_E= min\{\deg {i_P}_*\sigma^*E(\frak g/\frak p):(P,\sigma, \theta_P,U_P)$

is a Higgs reduction $\}$
\end{center}
 where, $i_P:U_P$ \ra $X $ is the inclusion. Since ${i_P}_*\sigma^*E(\frak
g/\frak p)$ is a
quotient of a fixed bundle $E(\frak g)$ on $X$, the integer $d_E$ is 
well defined.
Here we use the fact that the degrees of quotients of a fixed 
torsion-free
sheaf
is bounded from  below.                        
{\prop Let $(P,\sigma,\theta_P,U_P)$ be a Higgs reduction of structure
group
of $E$ to $P$ such that the following conditions hold: 
\begin{enumerate}
\item{} $\deg( \sigma^*E(\frak g/\frak p))=d_E$
\item{} $P$ is a maximal among all parabolic subgroups $P'$ satisfying
the condition that there is a Higgs reduction $(P', \sigma',\theta_{P'})$
of $E$ to $P'$ such that degree$(\sigma'^{*}E(\frak g/\frak p))=d_E$.
\end{enumerate}
Then the Higgs reduction $(P,\sigma,\theta_P,U_P)$ is a Higgs H-N
reduction.}

\pf Let  $(P,\sigma,\theta_P,U_P)$ be a Higgs reduction of $E$ to $P$
satisfying 
the properties stated in the proposition, $i.e$ we have a Higgs reduction
$(P,\sigma,\theta_P)$ of $E$ to
$P$ on the open subset $U_P$ of $X$. Let $U$ be the unipotent part of $P$. 
First note that $E_L$ as a $L$ bundle on $U_P$ gets a Higgs structure
induced by the surjective map $P \longrightarrow L \longrightarrow 0.$ 
We have, $\rho: E_P(\frak{p})\longrightarrow E_P(\frak{l})$ induced by the
above map, where ${\frak{p}}$ and ${\frak{l}}$ are Lie algebras of $P$ and
$L$ respectively.
So we define $\theta_L :=(\rho\otimes{id})\circ\theta_P$
where $\rho\otimes{id}
: E_P(\frak{p})\otimes\Omega^1_{U_P}\longrightarrow
E_P(\frak{l})\otimes\Omega^1_{U_P}$
and $\theta_P:\mathcal{O}_{U_P}\longrightarrow
E_P(\frak{p})\otimes\Omega^1_{U_P}$.\\      
 We first show that the associated Levi bundle $E_L$ is Higgs
semistable on $U_P$.\\
Suppose that $E_L$ is not a Higgs semistable $L$ bundle on
$U_P$. Therefore, there is a Higgs reduction  $(Q,\tau_Q,\theta_{Q},U')$
of $E_L$ to a parabolic subgroup $Q$ of $L$  such that,
\beqal{deg}
\deg(\tau^*T_{E_L/Q}) = \deg \tau^*E_L({\frak l/\frak q}) < 0
\eeqa                                                 
where, $U'$ is the open subset of $X$ with $ codim(U_P \backslash
U') \ge 2$ which implies that codim$(X \backslash U') \ge 2$. It is easy
to see that the inverse image of $Q$ under the projection
 $ P \longrightarrow L \longrightarrow 0 $  which will be denoted as
$P_1$, is a parabolic subgroup of $G$.
(Since, $G/P_1 \longrightarrow G/P$ is a fiber bundle with fiber $P/P_1 
\cong L/Q$, hence $P/P_1$ is complete. This implies that $G/P_1$ is
complete.)
Since $(Q,\tau_Q)$ is a reduction of structure group of $E_L$ to $Q$
given by the section $\tau_Q : U'\longrightarrow E_P(L/Q)$ we have a
section $\sigma_1: U'\longrightarrow E_P(P/P_1)$.  So we have
$(P_1,\sigma_1)$
a reduction of structure group of $E_P$ to $P_1$ on $U'$ open subset of
$U_P$, of $X$ with codim$(X \backslash U')\ge 2$. From now onwards we fix
this open subset $U'$ of $X$. Note that we have two short
exact sequences of $P$ and $P_1$ modules respectively.
\beqal{lie1} 
0\longrightarrow \frak{u} \longrightarrow \frak{p}\longrightarrow
\frak{l}\longrightarrow 0
\eeqa
\beqal{lie2} 
0\longrightarrow \frak u \longrightarrow\frak {p}_1\longrightarrow
\frak{q} \longrightarrow 0  
\eeqa
where $ \frak{u},\frak{p},\frak p_1,\frak{q}$ are Lie
algebras
of $U$,$P$,${P_1}$ and $Q$ respectively.
On $U'$ we have, $0 \longrightarrow E_{P_1}(\frak{p}_1)
\longrightarrow  E_P(\frak{p})$ and $0 \longrightarrow
E_{P_1}(\frak{q})
\longrightarrow  E_P(\frak{l})$. But we know that $ E_{P}(\frak{l})
=  E_L(\frak{l})$ and $ E_{P}(\frak{u})=  E_{P_1}(\frak{u})$ 
The above exact sequences therefore gives the following commutative diagram
of exact sequences.

\begin{center}
\begin{tabular}{ccccccccc}
&&$0$&&$0$&&$0$ \\
&& $\downarrow$ && $\downarrow$ && $\downarrow$ \\
0&$\longrightarrow$ & $E_{P_1}(\frak u)$ & $\longrightarrow $ &
$ E_{P_1}(\frak p_1) $ & $\longrightarrow$ & $E_{P_1}(\frak q)$ &
$\longrightarrow$ & 0 \\
&& $\downarrow$ && $\downarrow$ && $\downarrow$ \\
0& $\longrightarrow$ & $E_{P}(\frak u)$ & $\longrightarrow $ &
$E_{P}(\frak p) $ & $\longrightarrow$
& $E_{P}(\frak l)$ & $\longrightarrow $ & 0\\
\end{tabular}
\end{center}

This results in the following commutative diagram
of exact sequences.
\begin{center}
{\tiny 
\begin{tabular}{ccccccccc}
0&\ra&\hpauk&\ra&\hpapak&\rag&\hpaqk&\raga&\hapauk
\\
&&\uda&&\dai&&\daj&&\uda \\
0&\ra&\hpuk&\ra&\hppk&\raf&\hplk&\rafa&\hapuk \\
\end{tabular} 
} 
\end{center}
From the above diagram it follows that  there exists  $\theta_{P_1} \in
\mathrm {H}^0(U',E_{P_1}({\frak p}_1) \otimes \Omega^1_{U'})$ 
such that $i(\theta_{P_1})=\theta_{P}$.

Therefore we have got a Higgs structure $(P_1,\sigma_{1},\theta_{P_1})$
induced by 
Higgs structure $(Q,\tau_{Q},\theta_{Q})$ on $U'$.

But by the above diagram it is clear that this structure is compatible
with the Higgs structure induced by  $(P,\sigma ,\theta_P)$ 
i.e, we have the commutative diagram; 
\begin{center}
\begin{tabular}{ccc}
${\mathcal{O}}_{U'}$&$\stackrel{\theta_{P}}{\longrightarrow}$&
$E_{P}({\frak p}) \otimes \Omega^1_{U'}$ \\
&$\stackrel{\theta_{P_1}}{\searrow}$&$\uparrow$ \\
&&$E_{P_1}({\frak p}_1) \otimes \Omega^1_{U'}$ \\
&&$\uparrow$ \\
&&$0$  
\end{tabular} \\  
\end{center}
Let $\frak {p}_1,\frak {p},\frak {l}$ denote the Lie
algebras of $P_1,P,L$ respectively. We have the following exact sequence
of $P_1$ modules
\beqal{lie3}
0\longrightarrow \frak {p}/ \frak {p}_1\longrightarrow
\frak {g}/\frak {p}_1 \longrightarrow \frak {g}/\frak {p}
\longrightarrow 0   
\eeqa
Since$P/P_1 \cong L/Q$ 
we get $\frak {p} / \frak {p}_1 \cong \frak {l} / \frak {q}$ 
as $P_1$ modules.

Since $E_{P_1}$ is a Higgs reduction of $E$ to $P_1$, we consider the
vector bundles associated to $E_{P_1}$ using the $P_1$ modules in 
(\ref{lie3}). We have the following exact sequence of vector bundles 
on $U'$.
\beqa
0\longrightarrow \tau^{*} T_{{E_L}/Q}\longrightarrow \sigma^{*}
T_{{E}/P_1} \longrightarrow \sigma^{*} T_{{E}/P} \longrightarrow 0
\eeqa
Using (\ref{deg}) we conclude that
\beqa
\deg(\sigma_1^*T_{E/P_1})=\deg(\sigma^*T_{E/P})+ 
\deg(\tau^*T_{{E_L}/Q}) < \deg(\sigma^*T_{E/P}) = d_E
\eeqa
This contradicts the  assumption on $d_E$. Therefore, $E_L$ must be Higgs 
semistable.

Now we need to check the second condition in the definition of the  Higgs 
H-N reduction. Let $B$ be a Borel subgroup of $G$ contained in $P$ and $T
\subset B$ be the 
maximal torus. Let $\Delta$ be the system of simple roots. Let $I$ denote
the set of simple
roots defining the parabolic subgroup $P$.
So $I$ defines the roots of the Levi factor $L$ of $P$.
Take any dominant character $\chi$ of $P$ whose restriction to $T$ is
expressed as

\beqa
\chi\mid_T = \Sum_{\alpha\in\Delta-I}c_\alpha\alpha
+\Sum_{\beta\in I}c_\beta\beta
\eeqa

with $c_\alpha ,c_\beta \ge 0$.

Let $Z^0(L) \subset T$ be the connected component of the center of
$L$. The characters $\beta \in I$ of $T$ have the property that $\beta 
\mid_{{Z^0}(L)}$ is trivial. Hence we see that if $\chi$ is a
nontrivial character of $P$ of the above form then $c_\alpha > 0$ for some 
$\alpha\in\Delta-I$.

Also, some positive multiple of a character of $Z^0(L)$
extends to character of $L$.
Hence there are positive integers $n_\alpha$ for each $\alpha \in \Delta
-I$ such that $n_\alpha\alpha \mid_ {Z^0(L)}$ extends to a character
$\chi_\alpha$ on $L$. Now we see that $\chi_\alpha \mid_T$ can be
written as $n_\alpha \alpha + \sum_{\beta \in I} n_{\beta, \alpha} \beta$
for some integers  $n_{\beta, \alpha}$.

Let $N=\prod_{\alpha \in \Delta-I}n_\alpha$, then we
have the following equality.

\beqa
N\chi \mid_ T = \sum _{\alpha \in \Delta - I} Nc_\alpha\alpha 
+\sum _{\beta \in I} Nc_\beta\beta =  \sum _{\alpha \in \Delta - I}
c_\alpha '\chi _\alpha \mid_ T + \sum _{\beta \in I} c_\beta'\beta
\eeqa   

for some integer $c_\beta', \beta \in I$ and $c_\alpha'$ where $\alpha \in
\Delta - I$ such that $c_\alpha'$ is a positive multiple of $c_\alpha$.
Hence $N\chi-\sum_{\alpha\in\Delta - I} c_{\alpha}'\chi_{\alpha}$ is 
a 
character of $L$ whose restriction to $Z^0(L)$ is trivial.
Thus,$N\chi = \sum_{\alpha\in\Delta-I} c_\alpha'\chi_\alpha$.
Hence it is enough to prove the second condition for the characters 
of the form $\chi_\alpha$,
with $\alpha\in\Delta - I$. 

Fix an element $\alpha\in\Delta -I$. Let $P_2\supset P$ be the parabolic
subgroup of $G$ defined by the subset $I_2=\{\alpha\}\cup I$ of $\Delta$. 
Let $L_2$ be its Levi factor and $P'$ be the maximal parabolic subgroup 
of $L_2$ defined by the image of $P$ in $L_2$. Consider the group of all 
characters of $P'$ which are trivial on the center of $L_2$. This 
group 
is generated by a single dominant character $\omega$ of $P'$ with respect
to 
the root system of $L_2$ defined by its maximal torus $T$.

Now $\chi_\alpha$ defines a character on $P'$ which is trivial on center 
of $L$. Hence the restriction of $\chi_\alpha$ to $T$ can be written as 
$\chi_\alpha\mid_{T} =m\omega\mid_{T}$ for some integer $m$.

This enables us to write 
\beqa
m\omega \mid_{T} = n_\alpha \alpha + \sum_{\beta \in I} n_{\beta, \alpha}
\beta \nonumber.
\eeqa 
A dominant weight is always a nonnegative rational linear combination of 
simple roots. The fact that $n_\alpha > 0$ implies that $m > 0$ and 
$n_{\beta,\alpha} \geq 0$. Hence $\chi_\alpha$ is a positive multiple of 
$\omega$.

To verify the second condition of Higgs H-N reduction, we consider 
the surjective map $P_2\longrightarrow {L_2} \longrightarrow 0$ and the
injective map $P \stackrel {i} \hookrightarrow {P_2}.$ We note that
$P_2/P \cong L_{2}/P'$.
 Let $E_{P_2}$ be an extension of structure group of $E_P$ to
$P_2$, that is 
$E_P(P_2) = E_{P_2}$ on $U_P$. So $E_P$ is a reduction of structure
group of 
$E_{P_2}$ to $P$ on $U_P$. It is easy to see that $E_{P}(G) = E_{P_2}(G)$
and 
${P_2} \hookrightarrow G$ gives rise to a reduction  
$({P_2},\sigma_2)$of $E$ on $U_P$. Also we have $E_P(\frak
p)\stackrel{\rho'}
\hookrightarrow E_{P_2}({\frak p}_2)= (E_P({\frak p}_2))$ induced
by $i$.  

Now we define $\theta_{P_2}\colon=(\rho'\otimes
id)\circ\theta_P$

where $(\rho'\otimes id)\colon E_P(\frak p )\otimes\Omega^1_{U_P} 
\hookrightarrow E_P({\frak p}_2)\otimes \Omega^1_{U_P}$.

Since the following diagram commutes.
\begin{center}
\begin{tabular}{ccc}
${\mathcal{O}}_U$&$\stackrel{\theta_{P_2}}{\longrightarrow}$&
$E_{P}({\frak p}_2) \otimes \Omega^1_{U_P}$ \\
&$\stackrel{\theta_{P}}{\searrow}$&$\uparrow$ \\
&&$E_{P}({\frak p}) \otimes \Omega^1_{U_P}$ \\
&&$\uparrow$ \\
&&$0$   
\end{tabular} 
\end{center}  
we have a Higgs structure $\theta_{P_2}$ on  $E_{P_2}$
and hence Higgs reduction $(E_{P_2},\sigma_2,\theta_{P_2})$ on $U_P$.
As in the first part of the proposition $E_{L_2}=E_{P_2}(L_2)$ is a Higgs 
$L_2$ bundle.
Using the following exact sequence of $P$ modules,
\beqa
0 \longrightarrow\frak {p}_2/\frak p\longrightarrow\frak
g/\frak p\longrightarrow\frak g/\frak {p}_2\longrightarrow 0                                                           
\eeqa
we have the following exact sequence of vector bundles on $U_P$. 
\beqa
0\longrightarrow \sigma'^{*} T_{{E_{L_2}}/P'}\longrightarrow \sigma^{*}
T_{{E}/P} \longrightarrow {\sigma_2}^{*} T_{{E}/P_2} \longrightarrow 0
\eeqa 
where $\sigma'$ is the Higgs reduction of the structure group of
$E_{L_2}$ to $P'$ which is $E_P(P')$ for the obvious projection of $P$
to $P'$.
Here $E_P(P')$ gets Higgs structure in the following way.

 Let $f : E_P(\frak p) \longrightarrow E_P(\frak
p')$ be induced by the projection $P\rightarrow P'$. Define 
$\theta_{P'}\colon=(f\otimes id)\circ\theta_{P}$.
The section $\theta_{P'}$ is the required Higgs structure on $E_P(P')$.
Also note that the following diagram commutes.
\begin{center}
\begin{tabular}{ccc}
${\mathcal{O}}_{U_P}$&$\stackrel{\theta_{L_2}}{\longrightarrow}$&
$E_{L_2}({\frak l}_2) \otimes \Omega^1_{U_P}$ \\
&$\stackrel{\theta_{P}'}{\searrow}$&$\uparrow$ \\
&&$E_{P}({\frak p'}) \otimes \Omega^1_{U_P}$ \\
&&$\uparrow$ \\
&&$0$
\end{tabular} 
\end{center} 
From the above exact sequence we have 
\beqa
\deg(\sigma'^{*} T_{{E_{L_2}}/P'})=\deg(\sigma^{*}T_{{E}/P})-
\deg({\sigma_2}^{*} T_{{E}/P_2}). 
\eeqa
Now the assumption that $\deg(\sigma^{*}T_{{E}/P})=d_E$ gives the
inequality  $\deg(\sigma'^{*} T_{{E_{L_2}}/P'})<0$. Note that 
$det(\sigma'^{*}
T_{{E_{L_2}}/P'})$ is the line bundle associated to
the $P'$-bundle  $E_P(\mathrm P')$, for the character of $P'$ which is a
negative
multiple of $\omega$. Hence it follows from the above observation that
some 
positive powers of $L_{\chi_{\alpha}}^*$   and
$\det(\sigma'^{*}T_{{E_{L_2}}/P'})$ coincide, proving that
$\deg L_{\chi_{\alpha}} > 0$. Therefore, the second condition of the
definition of the canonical Higgs
reduction holds. This completes the proof of the proposition .
This proposition establishes the existence of a Higgs  canonical
reduction. 

\noindent
\section{ The Higgs structure on the adjoint bundle $E_G(\frak g)$}

Let $(E,\theta)$ be a principal Higgs $G$-bundle on $X$, let
$E_G(\frak g)$ denote 
its associated vector bundle for the adjoint representation  
$Ad: G \longrightarrow Gl(\frak g)$. That gives rise to a morphism of
Lie algebras $ ad: \frak g \longrightarrow End(\frak g)
 = \Lg \otimes\Lg^*$. So we get a morphism of vector bundles $\phi:\aeg
\longrightarrow \aeg \otimes \aeg^*$. We define $\adt =(\phi \otimes
id) \circ \theta: \mathcal {O}_X 
\longrightarrow  \aeg \otimes \aeg^* \otimes \Omega^1_X$ which gives Higgs
structure on $\aeg$, note that $\adt$ induces a homomorphism $\theta':\aeg
\longrightarrow \aeg \otimes \Omega^1_X$ of vector bundles on $X$.
Here we prove the following proposition  only for semiharmonic Higgs
principal $G$ bundles i.e. for Higgs principal bundles whose all chern
classes are vanishing. Therefore, the uniqueness of the H-N reduction 
holds only for semiharmonic  Higgs principal bundles which we prove in the
following section (c.f.\S 5) 

{\prop
\label {AD}
 A Higgs principal $G$ bundle $(E,\theta)$ on $X$  is Higgs
semistable iff the 
associated Higgs vector bundle $(\aeg,\adt)$ is Higgs semistable.}

\pf The proof of this
proposition goes along the same lines as that of \cite{AB}. 
First assuming that $(\aeg,\adt)$ is a Higgs semistable vector 
bundle, we need to prove that $(E,\theta)$ is a Higgs semistable principal $G$ bundle.
Suppose not, then there exists an open subset $U$, with
$\cod$ and a maximal parabolic subgroup $P$ of $G$ and a 
Higgs reduction $(P,\sigma,\theta_P,U)$ such that 
$\deg(\sigma^*(T_{E/P})) < 0$ (by Lemma \ref{para}), where $T_{E/P}$
is the relative tangent bundle for the projection $E \longrightarrow
E/P$. Let $E_0 =\aep \subset \aeg$ be the subbundle 
 given by the adjoint bundle of the $P$ bundle $E_P$ given by the 
reduction $\sigma$ then $\T=\aeg / E_0$. Also
($\aep,ad(\theta_P)$) is a Higgs subbundle of $(\aeg,\adt)$.
 Since $G$ is reductive, using a bilinear 
form which is non degenerate on $\frak g$, we have the identification 
$\aeg = \aeg^*$ which implies that $\deg(\aeg) = 0$.

 From the exact sequence of vector bundles on $U$
\beq
0\longrightarrow E_0\longrightarrow \aeg \longrightarrow \T \longrightarrow 0
\eeq
we conclude, $\deg(\aeg) = \deg(E_0) + \deg(\T)$. This implies that 
$\deg (E_0) > 0$ (since $\deg(\T) < 0$, $\mu(E_0) > 
\mu(E_G(\frak g))$ which contradicts the fact that ($\aeg,\adt$) is 
Higgs semistable.

Now we assume that $(E,\theta)$  is Higgs semistable then we will prove 
that $(\aeg,\adt)$ is a Higgs semistable vector bundle.

Suppose that $\aeg$ is not a Higgs semistable vector bundle on
$X$ and let
\beq
0=E_0 \subset E_1 \subset E_2\cdots\subset E_{k-1} \subset E_k=\aeg
\eeq
be the Higgs H-N filtration of the  Higgs vector bundle $\aeg$
(lemma 3.1. \cite{S2}).

For the sake of notational convenience, we denote $\aeg$ by $V$ and we
have the Higgs structure $\theta'$ induced by $\adt$ on $V$. Note that
we can choose an open subset $U$ of $X$ such that it contains all points
of codimension one (since $\cod)$ and every sheaf $E_i$ is locally 
free on $U$. We
fix this $U$. For any $x \in U$ we consider 
$E_{j,x}^\perp = \{v  \in V_x |< v,E_{j,x}> = 0$\}, where 
$<,>$ is the
non-degenerate bilinear form invariant under $G$ on  $V_x \cong \Lg$. 
Let $E_{j}^\perp$ be the kernel of the surjection $V \longrightarrow 
V^*$ (defined by the form) followed by canonical map $V^* 
\longrightarrow
 {E_j}^*$. 

We $claim$ that 
$E_{j}^\perp \cong (V/E_j)^*$ over $U$.  Let $v \in E_{j,x}^\perp$, then 
define $f: V_x \longrightarrow k$ by $f(v') = <v,v'>$. Since
$f(w) = <v,w>=0$ for all $w \in  E_{j,x}$. So $f$ induces a 
map 
$\bar{f}: V_x/E_{j,x} \longrightarrow k$. Thus $E_{j,x}^\perp
\hookrightarrow (V_x/E_{j,x})^*$, but dim$(E_{j,x}^\perp)$ =
dim$(V_x/E_{j,x})^*$, therefore $E_{j}^\perp \cong (V/E_j)^*$. 

We define $W_j: = E_{k-j}^\perp=(V/E_{k-j})^*$. Since $\theta' \mid
_{E_j}$ maps $ E_j \longrightarrow E_j \otimes \Omega^1_X$, we can define
$\bar{\theta'_{j}}:(V/E_j) 
\longrightarrow (V/E_j) \otimes  \Omega^1_X$ 
induced by $\theta': V \longrightarrow V \otimes\Omega^1_X$. 
So $((V/E_j)^*,\bar{\theta'_{j}})$ is a \emph{Higgs subbundle} of
$(V^*,\theta^1)$
where $\theta^1:\aeg^*\longrightarrow \aeg^*\otimes  \Omega^1_X$ induced
by $\adt$(since both $\theta'$ and
$\theta^1$ induced by $\adt$ and $(V/E_j)^*$ have Higgs structures 
induced
by $\theta'$). \\
Now we $claim$ that $W_i/W_{i-1}$ is a Higgs semistable
vector bundle and $\mu(W_i/W_{i-1})$ decreases as $i$
increases. To prove this claim, observe that $(E,\theta)$ is Higgs 
semistable vector 
bundle
if and only if 
$(E^*,\theta)$ is Higgs
semistable vector bundle.
Now we observe that,
\beqa
 W_{i}/W_{i-1} = 
\frac{(V/E_{k-i})^*}{(V/E_{k-(i-1)})^*}=(E_{k-(i-1)}/E_{k-i})^*
\eeqa 
 For, we have the diagrams:
\begin{center}
\begin{tabular}{ccccccccc}
&&&&&& 0&& 
\\
&&&&&& $\downarrow$ && \\
&&&&&& $E_{k-(i-1)}/E_{k-i}$&& \\
&&&&&& $\downarrow$ && \\
0&$\longrightarrow$&$E_{k-i}$ &$\longrightarrow$ &$V$ 
&$\longrightarrow$
&$ V/E_{k-i}$&$ \longrightarrow$&$ 0$ \\
&&&&$\Vert$ &&$\downarrow$& \\ 
0&$\longrightarrow $ &$ E_{k-(i-1)}$ &$ \longrightarrow$&$ V$ 
&$\longrightarrow$
&$ V/E_{k-(i-1)}$&$
\longrightarrow$&$ 0$ \\ 
&&&&&& $\downarrow$ && \\ 
&&&&&& 0&& 
\end{tabular}
\end{center}
\beqa
0\Ll (E_{k-(i-1)}/E_{k-i})^* \Ll (V/E_{k-i})^* \Ll  (V/E_{k-(i-1)})^* 
\Ll 0
\eeqa
and hence, $\frac{(V/E_{k-i})^*}{(V/E_{k-(i-1)})^*}$ =  
$(E_{k-(i-1)}/E_{k-i})^*$. Thus,
$W_{i}/W_{i-1} = (E_{k-(i-1)}/E_{k-i})^*$. 

But we know that $(E_{k-(i-1)}/E_{k-i})$ is Higgs semistable. 
Therefore $W_{i}/W_{i-1}$ is Higgs semistable. 
Note that 
$\mu(W_{i}/W_{i-1})=-\mu(E_{k-(i-1)}/E_{k-i})$, $\mu(W_{i+1}/W_{i})=
-\mu(E_{k-i}/E_{k-(i+1)})$. Hence, $\mu(W_{i}/W_{i-1})$ decreases as 
$i$ increases. This completes the proof of the claim. 

Therefore we obtain a Higgs H-N filtration of $V^*$ over $U$.
\beq
0=W_0 \subset W_1 \subset W_2 ......\subset W_{k-1} \subset W_k=V^*
\eeq 
such that each $W_i$ is a Higgs subbundle of $V$ and the quotient is Higgs
semistable; further, $\mu(W_{i}/W_{i-1})$ decreases as $i$ increases.

Now since $V \cong V^*$, we conclude that the filtration of $V^*$ over $U$
by the $W_{j}$'s coincides with the Higgs H-N filtration of $V$
(by the uniqueness of filtration). In other words, we have 
$E_{j} = E_{k-j}^\perp$ on $U$ for all $0 \leq j \leq k$.
Therefore, the above filtration is of the following form:
\beqa
 0 = E_{-l-1} \subset E_{-l} \subset E_{-l+1} \subset ....\subset 
E_{-1} 
\subset E_{0} \subset E_{1} \subset ...E_{l-1} \subset E_{l}=V
\eeqa
where $E_j$ is the orthogonal complement of $E_{-j-1}$ for the 
$G$-invariant form 
<,>.

 Let $\phi: E_{0} \otimes  E_{0} \lr V/E_{0},$ be the composition
of the Lie bracket operation with the natural projection
$V \lr E_{0}$. By the Proposition 2.9 in \cite{AB} we have, $\mu_{min} 
(E_0\otimes 
E_{0})=2 \mu_{min}(E_0)=2\mu (E_0/E_{-1})$
and $\mu_{max}(V/E_{0}) = \mu (E_1/E_0)$.

 From the fact that $E_{-1}$ is the orthogonal part of $E_0$, the
form <,> induces a non-degenerate quadratic form on 
$E_0/E_{-1}$. Consequently, we have $E_0/E_{-1} \cong (E_0/E_{-1})^*$ 
which implies that $\mu (E_0/E_{-1}) =0$.
 We have $\mu_{min}(E_0 \otimes E_0)=2\mu (E_0/E_{-1})=0 > 
\mu(E_1/E_0)=\mu_{max}(V/E_0)$. Therefore , it follows that
$\mathrm H^0(X,Hom_{Higgs}(E_0 \otimes E_0, V/E_{0})) =0$;
in particular, we have $\phi$ = 0. In other words, $E_{0} $ is closed 
under the Lie bracket. Consider
\begin{center} 
$\phi_{j} : E_{-j} \otimes E_{-1} \lr V/E_{-j-1}$ where $j \ge 0$ 
\end{center}
defined using the Lie bracket operation and the projection 
$V \lr V/E_{-j-1}.$

Repeating the above argument and using the property that 
$\mu (E_i/E_{i-1}) > \mu (E_{i+1}/E_{i})$ 
we deduce that $\phi_{j} =0$. In other words,
we have $[E_{-j},E_{-1}] \subset E_{-j-1}$ for any $ j \ge 0$.

 Using the above inclusion we conclude that $E_{-1}$ 
is a nilpotent Lie sub-algebra bundle of $E_0$, and $E_{-1,x}$ 
is also an ideal of $ E_{0,x}$ for any $x \in X$. 

Now by Lemma 2.11 of \cite{AB}, we see that over the open set 
$U$ of $X$, the sub-algebra bundle $E_0$ is a bundle of parabolic sub-algebras,
and it gives a reduction $\sigma:U\longrightarrow E_G/P$ of the structure
group of $E_G$ to a parabolic subgroup $P$ of $G$.

Using the above lemma we have a reduction $(P,\sigma)$. Let us denote the
principal $P$ bundle over $U$ obtained
in the above lemma by $E_P$  and $E_0 = E_P(\mr p)$. Since $E_0$ is a
Higgs subbundle of $\aeg$ we have,
$\theta'\mid_{E_0}:E_P(\frak p)\longrightarrow E_P(\frak p)\otimes
\Omega^1_U$ where $\theta' : E(\frak g) \rightarrow  E(\frak g) \otimes
\Omega^1_X$ induced by $\adt :\mathcal O_X \rightarrow \aeg \otimes \aeg^*
\otimes \Omega^1_X$. 


The section $\theta$  factors through $\aep \otimes \Omega^1_U$, call
it $\theta_P$. Therefore we conclude that $(E_P,\sigma_P,\theta_P)$ 
is a
\emph{Higgs reduction} of $E$ to $P$ on $U$.
Let $\chi_0$ be the character on $P$ associated to the action of
 $P$ on its Lie algebra $\mr p$. Then it is clear that $\chi_0$ is a
dominant character of $P$ with respect to the Borel subgroup $B$.
Hence from the Higgs semistability of $E_G$ we deduce that 
$\deg(\sigma^*E_G(\chi_0)) = \deg(E_P(\frak p)) = \deg(E_0) \le 0$ which
contradicts the fact that $\mu (E_0) > \mu (E) = 0$. This completes the
proof of the proposition.   

Now we will prove the  lemma that for every finite
dimensional linear representation of the group $G$ the associated bundle
is also Higgs semistable with some constraints. The proof of the lemma
goes along the same lines as that of the Lemma 5 of \cite{AAB} 
{\lem 
\label{Aut}A principal Higgs $G$-bundle $(E_G,\theta)$ over $X$ is Higgs
semistable
if and only if for every finite dimensional linear  representation 
$\rho:G \rightarrow Aut(V)$, such that $\rho (Z_0(G))$ is contained in the
center 
$Z(Aut(V))$, the associated Higgs vector bundle $E_G(V)$ is Higgs semistable.  
}

\pf Take the adjoint representation of $G$ on its
Lie algebra $\frak g$. Now by assumption the associated adjoint bundle is 
Higgs semistable, therefore by the previous proposition $(E_G,\theta)$ is
Higgs semistable.

Conversely, suppose that $(E_G,\theta)$ is Higgs semistable. By our
earlier discussion $E_G(V)$ gets Higgs section canonically; call it
$\theta_V$. Let $V=\bigoplus_{i=1}^{n}V_i$ be the decomposition of the 
$G$ module $V$ into irreducible submodules.\\
Consider the character $\bigwedge^{top}Hom(V_i,V_j)$ of $G$
given by the top 
exterior power. By assumption $Z_0(G)$ acts trivially on the line
$\bigwedge^{top}Hom(V_i,V_j)$. Since $G$ is connected, $G$ is a 
quotient(semi-direct product) of $Z_0(G) \times
[G,G]$, where $[G,G]$ denotes commutator subgroup of $G$. 
Any character of $[G,G]$
is trivial. So the action of $G$ on  $\bigwedge^{top}Hom(V_i,V_j)$ is 
trivial.

Let us denote the Higgs vector bundle $E_G(V_i)$ by
$(\mathcal W_i,\theta_i)$
associated 
to $E_G$ for the $G$ module $V_i$. So we have 
\beqa
\bigoplus_{i=1}^n \mathcal W_i = \mathcal W, &&   \bigoplus_{i=1}^n
\theta_i = \theta_V
\eeqa
By the earlier observation, the line bundle $\bigwedge^{top}Hom(\Wi,\Wj)$, 
which is associated to $E_G$ for the $G$ module  $\bigwedge^{top}Hom(V_i,V_j)$
is trivial. Therefore, we have $\mu(\Wi)= \mu(\Wj)$. This means that
$\mathcal W$ is Higgs
semistable iff each $\mathcal W_i$ is Higgs semistable.

So, our lemma is reduced to proving this only for irreducible
representations 
$\rho:G \longrightarrow Aut(V)$. Let $\mc W = E_G(V)$, we know that
$End(\mc W) = \mc W\otimes\mc{W}^*$ and
$\mc W = E_{GL(V)}(V)$ , then by Lemma \ref{G&V},
$\mc W$ is Higgs semistable iff  $E_{GL(V)}$ is Higgs semistable.
Again by proposition \ref{AD} $E_{GL(V)}$ is Higgs semistable iff
$ E_{GL(V)}(M_n(k)) = E_{GL(V)}
(V\otimes {V}^*) = \mc W\otimes\mc W^* = End(\mc W)$ is Higgs semistable.
Therefore, it is enough to prove that  $End(\mc W)$ is Higgs semistable.    
Since the $G$ module $V$ is irreducible, from Schur's lemma it follows
that the action of the center $Z(G)$ of $G$ on  $End(V)$ is trivial.
We note that the adjoint representation of $G$ on its Lie algebra
$\mr g$ gives a faithful representation of $G/Z(G)$. The group $G/Z(G)$
is also reductive. Therefore, the $G/Z(G)$ module $End(V)$ is a submodule
of the $G/Z(G)$ module  $\bigoplus_{j\in S}{\frak g}^{\otimes j}$, 
where $S$ is a finite collection of nonnegative integers
possibly with repetitions.

Therefore, the Higgs vector bundle $End(\mc W)$ is a Higgs subbundle of
$\bigoplus_{j\in S}({E_G(\frak g)})^{\otimes j}$.\\
Since $({E_G(\frak g)})$ is Higgs semistable of degree zero 
we have  ${(E_G(\frak g))}^{\otimes
j}$ is
Higgs semistable of degree zero
for every $j$. Thus, $\bigoplus_{j\in S}{(E_G(\frak g))}^{\otimes j}$
is
Higgs semistable of degree zero. But $End(\mc W)$ being a Higgs subbundle 
of degree zero, is also Higgs semistable. This completes the proof.
                                                              
{\lem Suppose $(E,\theta)$ is a Higgs principal $G$-bundle on $X$. Let
 $(E_P,\sigma_P,\theta_P,U_P)$ be a reduction, then $E_L$ is Higgs
semistable on $U_P$, where $L$ is a Levi factor of $P$, if and only 
if $i_*E_L(\frak l)$ is a Higgs
semistable torsion free sheaf over $X$.
}

{\pf Note that $E_L$  is Higgs semistable on $U_P$ if and only if
$E_L(\frak l)$  is Higgs semistable on $U_P$ (proposition \ref{AD}). 
But the
semistability of  $i_*E_L(\frak l)$ on $X$  is equivalent of the
semistability of $i_*E_L(\frak l)|_{U_P} = E_L(\frak l)$ on  $U_P$.

In conclusion, we have the following theorem.

{\thm Let $E_G$ be the Higgs semistable $G$ bundle. Let
$\rho: G\lr H$ be a
homomorphism of connected reductive groups
such that $\rho(Z_0(G))$ is contained in the connected component of the
center of $H$ containing the identity element. Let $E_H: =
E_G(H)$ be the Higgs principal $H$ bundle induced by  $\rho$. Then the $H$
bundle  $E_H$ is also Higgs semistable.}

{\pf Let $\phi: H\lr Aut(V)$ be a representation such that
$\phi(Z_0(H))\subseteq Z(Aut(V))$, where $Z_0(H)$ is the connected
component of the center of $H$ containing the identity.
By using Lemma \ref{Aut} it is enough to prove that
$E_H(V)$ is Higgs semistable. consider the composition $\phi\circ\rho:
G\lr Aut(V)$. Then $\phi\circ\rho(Z_0(G)\subseteq Z(Aut(V))$. Since $E_G$
is Higgs semistable, we have $E_G(V)$ is Higgs semistable.
But $E_G(V) = E_H(V)$. This proves the theorem.
 

\noindent
\section{Proof of the uniqueness of Higgs H-N reduction}
{\thm  Any Higgs $G$ bundle admits a unique H-N Higgs reduction.
In other words, if we fix a Borel subgroup $B$ of $G$, then there is an
unique H-N Higgs reduction to a parabolic subgroup $P$
containing $B$. 
}

\pf The existence of H-N Higgs reduction is proved in \S 2.
So, it is enough to prove the uniqueness.

 Let $(P_1,\sigma_{P_1},\theta_{P_1},U_{P_1})$ and
$(P_2,\sigma_{P_2},\theta_{P_2},U_{P_2})$
(where $P_1$ and $P_2$ are two parabolic subgroups containing $B$)
be two H-N Higgs reductions of the Higgs principal $G$ bundle
$(E,\theta)$.
 Note that we can choose open subset $U\subseteq X$ with $codim(X
\backslash U) \ge 2$ such that, the above two Higgs reductions exist on
this
$U$.
That is we have two H-N Higgs reductions of $E$ on $U$.
We will prove that the two Higgs subbundles $E_{P_1}(\frak p_1)$
and$E_{P_2}(\frak p_2)$ of $E_G(\frak g)$ are the same on $U$.
Equivalently we show $i_*(E_{P_1}(\frak p_1)) =
i_*(E_{P_2}(\frak p_2))$ on $X$, where $i$ is the inclusion of $U$ in $X$. 
And these are reflexive sheaves occurs 
in the H-N. filtration of $E(\frak g)$. For proving the above stated
result we need the following lemma.

{\lem Let $(E,\theta)$ be a principal $G$ bundle on a curve and
$(H_1,\sigma_{H_1}),
(H_2,\sigma_{H_2})$ be two reductions of $E$. 
Suppose the two subbundles $E_{H_1}(\frak h_1)$
and $E_{H_2}(\frak h_2)$ of $E(\frak g)$ coincide. Then there exists
an element $g\in G$ such that $H_1=gH_2g^{-1}$. Moreover, if the
normalizer of $H_1$ in $G$ is $H_1$ itself then the two reductions
$\sigma_{H_1}^*E$ and $\sigma_{H_2}^*E$ are the same.
}

The above lemma is easily seen to hold for nonsingular
projective varieties. By using this lemma we conclude that
$(P_1,\sigma_{P_1})=(P_2,\sigma_{P_2})$(since
$B$ is fixed and
normalizer of a parabolic subgroup is itself). Then by the definitions of
$\theta_{P_1}, \theta_{P_2}$
we arrive at the fact that
 $(P_1,\sigma_{P_1},\theta_{P_1}) = (P_2,\sigma_{P_2},\theta_{P_2})$.
So now our aim is to prove that the two Higgs subbundles
$E_{P_1}(\frak p_1)$ and $E_{P_2}(\frak p_2)$ of
 $E_G(\frak g)$
are the same.   
To prove this fact we need to prove one Lemma.

{\lem
\label{Hom}
Let $(V^1,\theta_1)$ and $(V^2,\theta_2)$ be two Higgs torsion-free
coherent sheaves on $X$ with Higgs filtrations
\beqa
0 &=& V_0\subset V_1\subset V_2\subset \cdots\cdots\cdots\subset
V_{k-1}\subset V_k=V^1\\                                                   
0 &=& V'_0\subset V'_1\subset V'_2\subset \cdots\cdots\cdots\subset
V'_{l-1}\subset V_l=V^2
\eeqa
where  $V_i/V_{i-1}$ and $V'_j/V'_{j-1}$  are Higgs semistable with 
 $\deg(V_i/V_{i-1}) \ge 0$ and $\deg(V'_j/V'_{j-1})< 0$  
for $0\le i\le k$ and $0\le j\le l$. Then $Hom_{Higgs}(V^1,V^2)=0$.
}

{\pf
Let $\phi : = V^1\longrightarrow V^2$ be a homomorphism of
Higgs torsion-free sheaves on $X$ with the above conditions.
We need to prove that $\phi=0$.
\noindent We choose an open set $U$ of $X$ with $codim( X \backslash U
) \ge 2
$
such that all the terms in the above two filtrations are locally free on
$U$. So it is enough to prove that $Hom_{Higgs}(V^1,V^2)=0$ of Higgs
vector bundles on $U$. To prove this fact we use the induction on the number of terms
in the filtrations i.e. on $k$ and $l$.

 Let $k=l=1$. Then $(V^1,\theta_1) , (V^2,\theta_2)$ are Higgs
semistable.
So we have, $ \mu (V^1) \le \mu(\phi (V^1)) \le \mu (V^2) $. 
But $\deg(V^1)\ge 0$ and $\deg(V^2)< 0$ by hypothesis,
so, $\mu(V^1) > \mu(V^2)$.
This leads to a contradiction.

We will prove the lemma for any $k$ and $l=1$.
Assume the lemma for $k'< k$ and $l=1$.  Note that we have 
\begin{center}
$\begin{array}{rccc}
0 \longrightarrow V_{k-1}\stackrel{i}{\longrightarrow}
&V_k&\stackrel{\pi}{\longrightarrow }
&V_k/V_{k-1} \longrightarrow 0\\
\phi_{|V_{k-1}}\searrow &\downarrow\,\phi& \swarrow \phi_1 \\
 &V^2&  &
\end{array} $        
\end{center}

So, $\phi \circ i = \phi\mid_{V_{k-1}} = 0 $. (by the induction on $k$)
Therefore we have $\phi_1: V_k/V_{k-1}\longrightarrow V^2 $
induced by $\phi$ such that
$\phi=\phi_1\circ \pi$. Again by induction
on $k$ we have $\phi_1=0$. This forces $\phi$ to be zero. 
So we have the lemma for $k$ and $l=1$.
Now assume the lemma for all $k$ and $l'< l$. Then we consider the
following,
\begin{center}
\begin {tabular}{rccc}
$0 \longrightarrow V'_{l-1}\stackrel{i'}{\longrightarrow}$&
$V'_l$&$\stackrel{\pi'}{\longrightarrow}$&$V'_l/V'_{l-1} \longrightarrow
0$\\
&$\uparrow\phi$ &$\nearrow \pi'\circ\phi$ \\
 &$V^1$&  &

\end {tabular}      
\end{center}
we get $\pi'\circ \phi=0$, that is im$(\phi)\subset V'_{l-1}$.

So, $\phi:V^1\longrightarrow V'_{l-1}$, this implies that $\phi =0$.
Thus, we have proved the lemma for all $k$ and all $l$.

Let $(P,\sigma_P,\theta_P,U_P)$ be a H-N Higgs reduction of
structure group of $E$ to $P$, we denote this bundle by $E_P$. The adjoint
action of $P$ on $\frak p$ preserves the sub-algebra $\frak u$.
The Higgs $L$ bundle $E_P(L)$, will be denoted by $E_L$.
By the definition of H-N Higgs reduction, $E_L$ is Higgs semistable.

The Higgs vector bundle associated to $E_P$ for the adjoint representation
of $P$ on $\frak p$ and $\frak g/\frak p$ are  naturally
identified with 
adjoint bundle of $E_P$ and $\sigma^*T_{E/P}$ respectively.

Consider the filtrations of $P$ modules
\beqa
0&=&V_0\subset V_1\subset V_2\subset \cdots\cdots\cdots\subset
V_{k-1}\subset V_k=\frak u\\
0&=&V'_0\subset V'_1\subset V'_2\subset \cdots\cdots\cdots\subset
V'_{m-1}\subset V'_m=\frak g/\frak p
\eeqa   
such that the quotients $W_i\colon=V_i/V_{i-1}$ and
$W'_j\colon=V'_j/V'_{j-1}$ are all  irreducible $P$ modules.

It follows that the action of $U$ on $W_i $( resp$W'_j)$ is trivial.
Therefore, the action of $P$ on $W_i ($resp$W'_j)$ factors through 
the quotient $L$. Let $\mathcal V_j, ($resp$\mathcal V'_j)$ denote 
the Higgs vector 
bundle over
$X$
associated to the $E_P$ for the $P$-module  $ V_j ($resp$ V'_j)$. Since $W_i$
and $W'_j$ are all irreducible $L$-modules and $E_L$ is Higgs 
semistable, $\mathcal W_i\colon=E_L(W_i)$ and $\mathcal
W'_i\colon=E_L(W'_i)$ are Higgs
semistable.
This result follows by Lemma \ref{Aut}. Let $B\subset P$ and $T$ a
maximal torus in $B$. Let $\Delta$ denote the
set of simple roots for $G$. Let $I\subset \Delta$ denote the set of simple
roots defining the parabolic subgroup $P$.
 The weights of $T$ on $\frak g/\frak p$ are of
the form 
\beqa
\gamma=\Sum_{\alpha\in\Delta}c_{\alpha}\alpha 
\eeqa
with $c_{\alpha}\le 0$ and  $c_{\alpha} < 0$ for at least one
$\alpha\in{\Delta - I}$. 
The weights of $T$ on $\frak u$ are of the form
$-\gamma$, where $\gamma$ is a weight on $\frak g/\frak p$.

From this it follows that the character of $P$ defined by the determinant
of the representation of $P$ on $W_j$ (resp, $W'_j$) is non-trivial and is
a non-negative (resp. nonpositive) linear combination of roots in 
$\Delta$.
Now by the second condition in the definition of a Higgs H-N
reduction we see that 
\beqa
\deg(\mathcal W_{j})>0,  && \deg(\mathcal W'_{j})<0
\eeqa

The adjoint action of $P$ on $\frak p/ \frak u$ factors through $L$
and this is precisely the adjoint representation of $L$. In other words,
the
Higgs vector bundle associated to $E_P$ for the $P$-module $\frak
p/ \frak u $ is $ E_L(\frak l)$. The Higgs semistability of $E_L$ implies
that
$E_L(\frak l)$ is Higgs semistable. Note that the degree of $E_L(\frak l)$
is zero.

Consider the exact sequence of Higgs vector bundles, 
\beqa
0\longrightarrow E_P(\frak p) \longrightarrow E_G(\frak g)\longrightarrow      
\sigma_P^*T_{E/P} \longrightarrow 0
\eeqa
corresponding to the following exact sequence of $P$-modules
\beqa
0\longrightarrow \frak p \longrightarrow \frak g\longrightarrow
\frak g/\frak p \longrightarrow 0
\eeqa
From the above observations it follows that the Higgs vector bundle
$\sigma_P^*T_{E/P}$ has the Higgs filtration of $\sigma_P^*T_{E/P}$ 
\beqal{tan1}
0=\mathcal V'_0\subset\mathcal V'_1\subset\mathcal V'_2\subset
\cdots\cdots\cdots\subset
\mathcal V'_{m-1}\subset \mathcal V'_m=\sigma_P^*T_{E/P}  
\eeqa
such that each quotient $\mathcal W'_i$ is Higgs semistable of negative
degree.

In the same way we have
\beqal{par1}
0=\mathcal V_0\subset\mathcal V_1\subset\mathcal V_2\subset
\cdots\cdots\cdots\subset
\mathcal V_{k-1}\subset \mathcal V_k \subset E_P(\frak p)
\eeqa
with each quotient  $\ W_i$ is Higgs semistable of positive degree
and $E_P(\frak p)/ \mathcal V_m$ which is identified with $E_L(\frak l)$,
Higgs semistable of degree zero.

So, we have the Higgs filtration of $E_P(\frak p)$ such that each quotient
is
of nonnegative degree.

Now we consider the exact sequences of Higgs vector bundles on $X$. 
\beqa
\label{Adjo1}
0\longrightarrow E_{P_1}(\frak p_1)\longrightarrow E_G(\frak 
g)\longrightarrow
\sigma_P^*T_{E/P_1}\longrightarrow 0
\eeqa
\beqa
\label{Adjo2}
0\longrightarrow E_P{_2}(\frak p_2)\longrightarrow
E_G(\frak g)\longrightarrow
\sigma_P^*T_{E/P_2}\longrightarrow 0
\eeqa
\noindent From (\ref{Adjo1}) and (\ref{Adjo2}) we have,
\beqa
E_{P_1}(\frak p_1)\longrightarrow\sigma^*T_{E/P_2} \\
E_{P_2}(\frak p_2)\longrightarrow\sigma^*T_{E/P_1}
\eeqa  
the homomorphisms of Higgs vector bundles.

Now by using Lemma \ref{Hom} 
we conclude that the above two homomorphisms are zero.
This implies that $E_{P_1}(\frak p_1)=E_{P_2}(\frak p_2)$.
So, finally we have proved the uniqueness of the H-N Higgs 
reduction.

\noindent
\section{H-N reduction for ramified G bundles over a smooth curve}
We assume from here onwards that $X$ is a smooth projective curve of genus
$g$, and $D=\sum_{i=1}^c D_i$ is a normal crossing reduced divisor. 
Unless stated otherwise, all the groups considered will be
reductive linear algebraic groups over $\mathbb{C}$ and $N \ge 2$ be a 
fixed integer.

{\defi A ramified $G$ bundle of type $N$ on $X$ with ramification over
$D$ is a smooth variety $Q$ with an action of $G$ and a morphism
$\phi: Q \lr X$ satisfying the following properties:
\begin{enumerate}
\item{} The action $G$ on $Q$ is proper.
\item{} The triple $(Q,\phi,X)$ is a geometric quotient.
\item{} In the complement of the divisor $D\subset X$ the morphism 
$\phi: Q \lr X$ is a principal $G$ bundle.
\item{} At the finitely many orbits on which the isotropy is nontrivial,
the isotropy is a cyclic subgroup whose order divides $N$.
\end{enumerate}
}
Given a homomorphism $\rho : G \lr H $ and a ramified $G$-bundle $Q$ of
type $N$ on $X$, the 
quotient space $Q\times_G H$ has a natural structure of a ramified
$H$-bundle of type $N$ on $X$. This construction is called the extension
of the structure group of $Q$ to $H$. Let us denote the quotient space by 
$Q(H)$(c.f. \cite {BBN1}).

{\defi Let $\phi_1 : Q \lr X$ and $\phi_2 :Q' \lr X$ be two ramified
$G$-bundles of type $N$ over $X$. Let $h : Q \lr Q'$ be a morphism of
varieties. We say  $h$ is a
morphism of ramified $G$-bundles if
\begin{enumerate}
\item{} $\phi_2\circ h = \phi_1$
\item{} The isotropy for the image of the orbits(whose isotropy is
nontrivial) is a cyclic subgroup whose order divides $N$ 
\end{enumerate}         
}             
We recall some of the definitions stated in \cite{BBN} and 
\cite{BBN1}.
\begin{enumerate}
\item{} Let $H$ be a closed subgroup of $G$. A reduction of structure
group
of a ramified $G$-bundle $Q$ to $H$ is given by a section $\sigma : X \lr 
Q/H$. Let $Q_H$ be the inverse image of the natural projection of $Q$ to
$Q/H$, of the subset of $Q/H$ defined by the image of $\sigma$. This $Q_H$
is a ramified $H$-bundle on $X$. It is easy to see that $Q_H \times_H G
\cong Q$.
\item{} Let $Q_H$ be a reduction of a ramified $G$-bundle $Q$. If $W$ 
is a finite dimensional $H$-module, then the associated construction 
$Q_H {\times}^H W$, is a parabolic vector bundle over $X$. This is 
denoted by $Q(W)_*$.
\item{} A ramified $G$-bundle $Q$ is called semistable if for every
reduction of structure group $(P,\sigma)$ of $Q$ to any parabolic subgroup
 $P$ of $G$ where $\sigma:X \lr Q/P$ is a section, and any dominant
character $\chi$ of $P$, the parabolic line bundle
$Q(\chi)_*$ $(=L_{\chi}$, as in our first section) has non-positive 
parabolic
degree.  
\item{} Let $Q$ be a ramified $G$ bundle of type $N$ over
$X$ with ramification divisor $D$. A Kawamata cover of type $(N,D)$
(or simply N) is a Galois cover $p : Y \stackrel{\Gamma}{\lr} X$
where $Y$ is a connected smooth projective variety and $\Gamma$ is a
finite subgroup of $Aut(Y)$. A $(\Gamma,G)$-bundle is defined to be a 
principal $G$-bundle $E$ together with a lift of the action 
of $\Gamma$ on $E$ which commutes with the right action of $G$ on 
$E$.
\item{}The ramified $G$-bundle $Q$ over $X$ with ramification divisor 
$D$ gives rise to a $(\Gamma,G)$-bundle $E$ over $Y$ such that for 
every $y \in (p^*D_i)_{red}$ has a cyclic subgroup $\Gamma_y \subset 
\Gamma$
of order $k_{i}N$ as the isotropy group.
\item{} Let [$\Gamma,G,N$] denote the collection of 
$(\Gamma,G)$-bundles on $Y$ satisfying the following two conditions: 
\begin{enumerate}
\item{} For a general point $y$ of an irreducible component of 
$(p^*D_{i})_{red}$, the action of $\Gamma_y$ on $E_y$(fibre over $y$)
is of order $N$.
\item{} For a general point $y$ of an irreducible component of a 
ramification divisor for $p$ not contained in $(p^*D)_{red}$, the 
action of $\Gamma_{y}$ on $E_y$ is trivial.
\end{enumerate}
\item{} Let $E \in [\Gamma,G,N]$ and $E'=E/\Gamma$. The variety $E'$ 
is smooth; 
moreover $\phi : E' \lr X $ is a ramified $G$-bundle.
\item{} $E$ is a normalisation of the fibre product $Y\times_X E'$. In
other words, the $G$-bundle $E$ on $Y$ can be constructed back from $E'$.
\item{} A ramified $G$-bundle $Q$ canonically defines a parabolic
$G$-functor $F_Q : Rep(G) \lr PVect( X, D)$ taking values in a category
$PVect( X, D, N)$ for some $N$. Conversely, given a parabolic $G$-functor
$F$ there exists a ramified $G$-bundle $Q$, unique up to isomorphism of
ramified $G$-bundles, such that $F_Q = F$. 
\item{} A ramified $G$-bundle $Q$ is semistable if and only if the
corresponding functor $F_Q$ is semistable and hence if and only if the
induced $(\Gamma,G)$-bundle is semistable. A ramified $G$-bundle $Q$ is
stable if and only if the induced $(\Gamma,G)$-bundle is stable.
\end{enumerate}  

{\defi Let $Q$ be a ramified $G$-bundle over $X$. Let $(P ,\sigma)$ be a 
reduction of structure group of $Q$ to a parabolic subgroup $P$ of $G$.
Then the reduction  $(P ,\sigma)$ is called a H-N reduction for the
ramified $G$-bundle $Q$ on $X$ if the following two conditions hold:
\begin{enumerate}
\item{} Let $L$ be the Levi factor of $P$. Then the ramified $L$-bundle
$Q_P(L)$ is semistable.
\item{} For any dominant character $\chi$ of $P$ with respect to some
Borel subgroup $B\subset P$ of $G$ the associated parabolic line bundle
$Q_P(\chi)_*$  has parabolic degree $> 0$.
\end{enumerate} 
}
\noindent
\begin{subsection} 
{H-N reduction for $[\Gamma,G,N]$ bundle} 
We define the H-N reduction for a
$[\Gamma,G,N]$ bundle over any connected smooth projective variety 
over $\mathbb{C}$(c.f.\cite{BBN}) and we prove that it exists and is
unique.
 
{\defi A $(\Gamma,G,N)$ bundle $E$ over $Y$ with a $\Gamma$-reduction 
of structure
group
$(P,\sigma)$ over some open subset $U$ with $\cod$
is said to be canonical or \emph {$\Gamma$ equivariant H-N reduction} if the
following two conditions hold:
 
\begin{enumerate}
\item{} The associated Levi bundle $E_P(L)$,
 is $\Gamma$-semistable.
 
\item{} For every dominant character $\chi$ of $P$ with respect to some
Borel subgroup contained in $P$, the line bundle $E_P(\chi)$ on $U_P$
has positive degree.
\end{enumerate}
}
 
{\prop For any ($\Gamma$,G)-bundle over $Y$, $\Gamma$-equivariant H-N 
reduction exists and it is unique with respect to some Borel subgroup 
$B$ contained in $P$.                                                                      
}
 
\pf We will use the fact that for  any $G$-bundle the canonical 
reduction
(P,$\sigma$) exists and the usual semistability is same as
$\Gamma$-semistability(\cite{BBN}, Proposition 4.1). So our main job is to
show that the usual H-N reduction (P,$\sigma$) is $\Gamma$-saturated.
Let $\gamma \in \Gamma (\subset Aut(Y))$, we note that $\gamma$ 
induces \Gt $:E \lr E$. 
We denote $\sigma^*E$ by 
$E_P$ and \Gt($E_P$) by $\Et$. Clearly $\Et$ is a
$P$ bundle. [We define P action ``$\star$'' on $\Et$ by $p
 \star$($\gamma$(e)):=$\gamma(e) \cdot p^{-1} = \gamma (e \cdot p^{-1})$ 
, where $\cdot$ denotes the action of $G$ on $E$,
clearly $\gamma(p \star e) =p \star (\gamma(e))$. This means the free 
action
of ``P'' on $\Et$ commutes with the action of $\gamma$. Hence $\Et$ 
is a 
$P$-bundle.] \\

Note that $E_P \cong \Et$, by uniqueness of H-N reduction
$E_P$ = $\Et$. This implies that $E_P$ is $\Gamma$
saturated (since $\gamma$ is arbitrary). So $E_P$ is a $(\Gamma,P)$ reduction.
Now we use the fact that that usual semistability of $\Et(L)$ is same
as $\Gamma$ semistability of $\Et(L)$(c.f. Theorem 4.3, \cite 
{BBN1}).
The second condition for $\Gamma$ equivariant H-N reduction is satisfied 
automatically.
 
So $E_P$ is a \emph{$\Gamma$ equivariant H-N                                          
reduction}, and it is unique.

\end{subsection}
\noindent
\begin{subsection} 
{Proof of H-N reduction for ramified $G$ bundle}

We fix a Borel subgroup $B$. We would like to prove here that 
for any ramified $G$-bundle there exists 
unique H-N reduction. We already know that the existence of the H-N
reduction for the $(\Gamma, G, N)$ bundles. Let $\phi : Q \lr X$ be a
ramified $G$-bundle on $X$ and let $p : Y \lr X$ be a
Kawamata cover and $E \in [\Gamma, G ,N]$ such that $E/\Gamma = Q$.
Let $(P,\sigma)$ be a H-N reduction of $E$ to a parabolic subgroup $P$ 
of $G$ where $\sigma : Y \lr E(G/P)$ is the $\Gamma$-equivariant section. 
We will prove that the corresponding reduction $\tau : X \lr Q/P$ of 
$Q$ to
$P$  (where $\tau$ is induced by the $\Gamma$-equivariant section $\sigma
$)is the H-N reduction for the ramified $G$-bundle $Q$.
Note that the semistablity of $E_P(L)$ is equivalent to the semistablity
of the ramified $L$-bundle $(E_P/\Gamma)(L) = E_P(L)/\Gamma$.
Let $\chi$ be the dominant character of $P$ with respect to the Borel
subgroup $B\subset P$. It is proved in \cite {BBN} that $\deg(E_P(\chi))
\ge$ $0$ if and only if $\deg((E/\Gamma)(\chi)_*) \ge 0$. Therefore, 
by the
above discussion it is clear that the
reduction $(E_P/\Gamma ,\tau)$ is a H-N reduction for the ramified
$G$-bundle $Q = E/\Gamma$ .

The uniqueness of the H-N reduction follows from the uniqueness of the 
H-N reduction of the $(\Gamma, G)$-bundle $E$.
\end{subsection}


\noindent



 \noindent

\begin{thebibliography}{99}
\bibitem{AAB} Anchouche.B,H.Azad,Biswas.I Harder Narasimhan Reduction for
a principal bundle over a compact K$\ddot{a}$hler manifold, Math.Ann.
 {\bf 323},693-712 (2002).
\bibitem{AB}  Anchouche.B, Biswas.I: Einstein-Hermitian
connections on
polystable principal
bundles over a compact K$\ddot{a}$hler manifold. Amer. Jour. Math. {\bf
123},
207-228 (2001).
\bibitem{AtB} Atiyah,M.F., Bott, R.:The Yang-Mills equations over Reimann 
surfaces. Phil.Trans.R.Soc. Lond.A {\bf 308}, 523-615 (1982)

\bibitem{BBN} Balaji.V, Biswas.I,Nagaraj.D.S: Ramified G-bundles as 
parabolic bundles, J. Ramanujan Math. Soc.  {\bf 18}, 123-138,No.2 (1993)
\bibitem{BB} Balaji.V, Biswas.I: On principal bundles with vanishing Chern
classes, J.Ramanujan Math.Soc. {\bf 17}, No.3 187-209(2002).
\bibitem{BBN1}  Balaji.V, Biswas.I, Nagaraj.D.S: Principal bundles over
projective manifolds with parabolic structure over a divisor,
Tohuku. Math.
Journal, {\bf 53}, 337-367, (2001)  
\bibitem{BH} Biswas.I, Holla.Y: Harder Narasimhan Reduction of a principal
bundle, preprint.
\bibitem{BR} Brhend.K.A : Semistability of reductive group schemes over 
curves; Math. Ann, {\bf 301}, 281-305 (1995)         
\bibitem{D} Deligne, Milne, Ogus, Shih; Hodge cycles, Motives, and Shimura
varieties; Springer-Verlag, LNM {\bf 900}.
\bibitem{F} Faltings.G, Stable G-bundles and projective connections,
J. Algebraic Geometry. {\bf 2}, 507-568, (1993). 
\bibitem{HN} Harder.G, Narasimhan.M.S: On the cohomology groups of moduli
spaces of vector bundle on curves, Math Ann. {\bf 212}, 215-248 (1975).  
\bibitem{H} Humphreys.J: Linear algebraic groups GTM 21, Springer verlag,
1975.                   
\bibitem{L}  Langton.S.G, Valuative criterion for families of vector
bundles on algebraic varieties, Annals of Mathematics {\bf 101} 88-110 
(1975).
\bibitem{MS} Mehta.V.B, Seshadri.C.S: Moduli of vector bundles on
  curves with parabolic structure, Math Ann. {\bf 248}, 205-239
  (1980).

\bibitem{AR1} Ramanathan.A: Stable principal bundles on a compact Riemann
surface. Math
Ann {\bf 213} 129-152 (1975).
\bibitem{AR2} Ramanathan.A: Moduli of principal bundles, Algebraic
geometry         
proceedings, Copenhagen LNM {\bf 732} 527-532 (1979). 
\bibitem{S1} Simpson.C : Higgs Bundles and local systems,
Publ. Inst. Hautes. Etudes. Sci. {\bf 75} 5-95 (1992).
\bibitem{S2} Simpson.C : Moduli of representations of the fundamental
group of a smooth projective
variety-1, Publ. Inst. Hautes. Etudes. Sci. {\bf 79} 47-129 (1994). 
\end{thebibliography}
\end{document}